\documentclass[journal]{IEEEtran}

\usepackage[utf8x]{inputenc}
\usepackage{amssymb,amsmath,amsbsy}
\usepackage{authblk}
\usepackage[usenames,dvipsnames]{color}
\usepackage{hhline}
\usepackage{float}
\usepackage{subfig}
\usepackage{algorithm}
\usepackage{algorithmic}
\usepackage{graphicx}	
\usepackage{array,multirow}
\usepackage{float}
\usepackage{soul,color} 
\usepackage{mathtools}
\usepackage{cite}
\usepackage[]{units} 	

\newcommand{\set}{ \leftarrow }

\newcommand{\expt}{{\mathbb E}}

\DeclareMathOperator*{\argmin}{arg\min}
\newcommand{\child}{\operatorname{child}}
\newcommand{\parent}{\operatorname{anc}}
\newcommand{\nodes}{\mu}
\renewcommand{\Re}{\mathbb{R}}
\newcommand{\N}{\mathbb{N}}
\newcommand{\prox}{\operatorname{prox}}
\newcommand{\proj}{\operatorname{proj}}
\newcommand{\dist}{\operatorname{dist}}

\renewcommand{\argmin}{\operatorname*{argmin}}

\newcommand{\dfn}{\mathrel{\mathop:}=}

\begin{document}

\title{GPU-accelerated stochastic predictive control of drinking water networks}

\author{Ajay~K.~Sampathirao, Pantelis~Sopasakis, Alberto~Bemporad~\IEEEmembership{Fellow,~IEEE} and Panagiotis~Patrinos
\thanks{The first three authors are with IMT Institute for Advanced Studies Lucca, 
Piazza San Francesco 19, 55100 Lucca, Italy.
Emails: \texttt{\{a.sampathirao, p.sopasakis, a.bemporad\}@imtlucca.it}.}
\thanks{P.~Patrinos is with STADIUS Center for Dynamical Systems, Signal Processing and Data Analytics,
KU Leuven, Department of Electrical Engineering (ESAT), Kasteelpark Arenberg 10, 
3001 Leuven, Belgium. Email: \texttt{panos.patrinos@esat.kuleuven.be}.}
\thanks{This paper was submitted on \texttt{arXiv} on 4 April 2016.}}


\maketitle

\begin{abstract}
Despite the proven advantages of scenario-based stochastic model predictive 
control for the operational control of water networks, its 
applicability is limited by its considerable computational footprint.
In this paper we fully exploit the structure of these problems and solve them
using a proximal gradient algorithm parallelizing the involved operations. 
The proposed methodology is applied and validated on a case study: the water network of 
the city of Barcelona.
\end{abstract}



\section{Introduction}

\subsection{Motivation}
Water utilities involve energy-intensive processes, complex in nature (dynamics) and
form (topology of the network), of rather large scale and with interconnected components,
subject to uncertain water demands from the consumers and are required to 
supply water uninterruptedly.
These challenges call for operational management technologies able to provide
reliable closed-loop behavior in presence of uncertainty.
In 2014, the IEEE Control Systems Society identified many aspects of 
the management of complex water networks as emerging future research directions~\cite{CST14}.

Stochastic model predictive control is an advanced control scheme which can 
address effectively the above challenges and has already been used for the management
of water networks~\cite{GroOcam+14,GroMaeOca14}. However, unless restrictive assumptions
are adopted regarding the form of the disturbances, such problems are known to 
be computationally intractable~\cite{GroMaeOca14,NemSha06}. In this paper we combine 
an accelerated dual proximal gradient algorithm with general-purpose graphics processing 
units (GPGPUs) to deliver a computationally feasible solution for the control of 
water networks.

\subsection{Background}
 
The \textit{pump scheduling problem} (PSP) is an optimal control problem for determining an open-loop 
control policy for the operation of a water network. Such open-loop approaches are known since the 
80's~\cite{Crea88,ZesSha89}. More elaborate schemes have been proposed such as~\cite{YuPoSt94} where
a nonlinear model is used along with a demand forecasting model to produce an optimal open-loop
24-hour-ahead policy. Recently, the problem was formulated as a mixed-integer nonlinear program to
account for the on/off operation of the pumps~\cite{BagBar+13}. Heuristic approaches using evolutionary algorithms, 
genetic algorithms, and simulated annealing have also appeared in the literature~\cite{McCPow04}.  
However, a common characteristic and shortcoming of these studies is that they assume 
to know the future water demand and they do not account for the various sources of uncertainty which 
may alter the expected smooth operation of the network.

The effect of uncertainty can be attenuated by feedback from the network
combined with the optimization of a performance index taking into account the 
system dynamics and constraints as in PSP. This, naturally, gives rise to
model predictive control (MPC) which has been successfully used 
for the control of drinking water networks~\cite{SamGroSop+14,OcPuCe+09}.
Recently, Bakker \textit{et al.} demonstrated experimentally on five full-scale 
water supply systems that MPC will lead to 
a more efficient water supply and better water quality than a conventional 
level controller~\cite{BaVrePa+13}. Distributed and decentralized MPC formulations
have been proposed for the control of large-scale water networks~\cite{LeZaNe+10,OcFaBa+10}
while MPC has also been shown to be able to address complex system dynamics such as
the Hazen-Williams pressure-drop model~\cite{SanKumNarNarBha15}.

Most MPC formulations either assume exact knowledge of the system
dynamics and future water demands~\cite{OcPuCe+09,OcFaBa+10} or endeavor to
accommodate the worst-case scenario~\cite{SamGroSop+14,TraBrd09,GorNEm14,WatMcK97}.
The former approach is likely to lead to adverse behavior in presence of
disturbances which inevitably act on the system while the latter turns 
out to be too conservative as we will later demonstrate in this paper.

When probabilistic information about the disturbances is available 
it can be used to refine the MPC problem formulation. The uncertainty
is reflected onto the cost function of the MPC problem deeming it 
a random variable; in \textit{stochastic MPC} (SMPC) the index to minimize is typically
the expectation of such a random cost function under the (uncertain)
system dynamics and state/input constraints~\cite{CaKoWu09,BerBem12}.

SMPC leads to the formulation of optimization problems over spaces of
random variables which are, typically, infinite-dimensional.
Assuming that disturbances follow a normal probability distribution
facilitates their solution~\cite{HesBos02,BerBro07,NemSha06}; 
however, such an assumption often fails to be realistic.
The normality assumption has also been used for the stochastic control of 
drinking water networks aiming at delivering high quality of services
-- in terms of demand satisfaction -- while minimizing the pumping cost
under uncertainty~\cite{GroOcam+14}.

An alternative approach, known as \textit{scenario-based
stochastic MPC}, treats the uncertain disturbances as discrete random variables
without any restriction on the shape of their distribution%
~\cite{CalCam06,CaGaPr09,PrGaLy12}.
The associated optimization problem in these cases becomes a
discrete multi-stage stochastic optimal control problem~\cite{ShaDen09}.
Scenario-based problems can be solved algorithmically, however, 
their size can be prohibitively large making them impractical 
for control applications of water networks as pointed out by Goryashko and Nemirovski~\cite{GorNEm14}.
This is demonstrated by Grosso \textit{et al.} who provide a comparison
of the two approaches~\cite{GroMaeOca14}.
Although compression methodologies have been proposed -- such as 
the scenario tree generation methodology of Heitsch and R{\"o}misch~\cite{HeiRom09} --
multi-stage stochastic optimal control problems may still involve 
up to millions of decision variables.

Graphics processing units (GPUs) have been used for the acceleration of the
algorithmic solution of various problems in signal processing~\cite{McCool07},
computer vision and pattern recognition~\cite{BenMatTim+12}
and machine learning~\cite{Jang+08,Guzhva2009} leading to a manifold 
increase in computational performance. To the best of the authors' knowledge,
this paper is the first work in which GPU technology is used for the solution of
a stochastic optimal control problem.

\subsection{Contributions}
In this paper we address this challenge by devising an optimization
algorithm which makes use of the problem structure and sparsity.
We exploit the structure of the problem, which is dictated by the 
structure of the scenario tree, to parallelize the involved operations.
Then, the algorithm runs on a GPU hardware leading to a significant
speed-up.

We first formulate a stochastic MPC problem using a linear 
hydraulic model of the water network while taking into
account the uncertainty which accompanies future water demands.
We propose an accelerated dual proximal gradient algorithm
for the solution of the optimal control problem and report results
in comparison with a CPU-based solver.

Finally, we study the performance of the closed-loop system
in terms of quality of service and process economics using 
the Barcelona drinking water
network as a case study. We show that the number of scenarios
allows us to refine our representation of uncertainty and trade
the economic operation of the network for reliability and quality 
of service.

\subsection{Mathematical preliminaries}
Let $\bar{\Re}=\Re\cup \{+\infty\}$ denote the set of 
extended-real numbers.
The set of of nonnegative integers $\{k_1,k_1+1,\ldots, k_2\}, k_2\geq k_1$ 
is denoted by $\mathbb{N}_{[k_1,k_2]}$. 
For $x\in\Re^n$ we define $[x]_+$ to be the vector in $\Re^n$ whose
$i$-th element is $\max\{0, x_i\}$. For a matrix $A\in\Re^{n\times  m}$
we denote its transpose by $A'$.

The \textit{indicator function} of a set $C\subseteq \Re^n$ 
is the extended-real valued function $\delta(\cdot|C):\Re^n\to\bar{\Re}$ and 
it is $\delta(x|C)=0$ for $x\in C$ and $\delta(x|C)=+\infty$ otherwise. 
A function $f:\Re^n\to \bar\Re$ is called \textit{proper} if there is a $x\in\Re^n$
so that $f(x)<\infty$ and $f(x)>-\infty$ for all $x\in\Re^n$.
A proper convex function $f:\Re^n\to \bar\Re$ is called \textit{lower semi-continuous} or 
\textit{closed} if for every $x\in\Re^n$, $f(x)=\liminf_{z\to x}f(z)$.
For a proper closed convex function $f:\Re^n\to \bar\Re$, we define its 
\textit{conjugate} as $f^*(y) = \sup_x \{ y'x - f(x)\}$.
We say that $f$ is \textit{$\sigma$-strongly convex} if $f(x)-\frac{\sigma}{2}\|x\|_2^2$ is 
a convex function. Unless otherwise stated, $\|\cdot\|$ stands for the Euclidean norm.

\section{Modeling of Drinking Water Networks}\label{S.2}

\subsection{Flow-based control-oriented model}
Dynamical models of drinking water networks
have been studied in depth in the last two 
decades~\cite{OcPuCe+09,MiFu84,OcFaBa+10}.
Flow-based models are derived from simple mass balance equations of the network
which lead to the following pair of equations
\begin{subequations}
\begin{align}
x_{k+1}&=A x_{k} + B u_{k} + G_d d_{k},\label{eqn:model}\\
0&=E u_k+E_d d_k,\label{eq:demand-control}
\end{align}
\end{subequations}
where $x \in \Re^{n_x}$ is the state vector corresponding to the volumes 
of water in the storage tanks, $u\in\Re^{n_u}$ is the vector
of manipulated inputs and $d\in\Re^{n_d}$ is the vector of water demands. 
Equation~\eqref{eqn:model} forms a linear time-invariant system with additive
uncertainty and~\eqref{eq:demand-control} is an algebraic input-disturbance 
coupling equation with $E\in\Re^{n_e\times n_u}$ and $E_d\in\Re^{n_e\times n_d}$ 
where $n_e$ is the number of \textit{junctions} in the network.

The maximum capacity of the tanks and the maximum pumping capacity of each pumping 
station is described by the following bounds:
\begin{subequations}
\begin{align}
u_{\min}\leq &u_{k} \leq u_{\max},\\
x_{\min}\leq &x_{k} \leq x_{\max}.
\end{align}
\end{subequations}

The above formulation has been widely used in the 
formulation of model predictive control
problems for DWNs~\cite{SamGroSop+14,GroOcam+14,OcFaBa+10}.

\subsection{Demand prediction model}
The water demand is the main source of uncertainty that affects
the dynamics of the network. Various time series models have been proposed 
for the forecasting of future water demands such as seasonal Holt-Winters, 
seasonal ARIMA, BATS and SVM~\cite{SamGroSop+14,Wang2014}. 
Such models can be used to predict nominal forecasts of the upcoming water
demand along a horizon of $N$ steps ahead using measurements available up 
to time $k$, denoted by $\hat{d}_{k+j\mid k}$. 
Then, the actual future demands $d_{k+j}$ --- which are unknown to the controller
at time $k$ --- can be expressed as
\begin{equation}\label{eq:demand:model}
d_{k+j}(\epsilon_j)=\hat{d}_{k+j\mid k}+\epsilon_j, 
\end{equation}
where $\epsilon_j$ is the demand prediction error which is a random 
variable on a probability space $(\Omega_{j},\mathfrak{F}_j, \mathrm{P}_j)$
and for convenience we define the tuple $\boldsymbol{\epsilon}_j=(\epsilon_0,\epsilon_1,\ldots,\epsilon_j)$
which is a random variable in the product probability space.
We also define $\hat{\mathbf{d}}_k=(\hat{d}_{k\mid k},\ldots\hat{d}_{k+N-1\mid k} )$.

\begin{figure}[ht]
 \centering
 \includegraphics[keepaspectratio=true,width=0.49\textwidth]{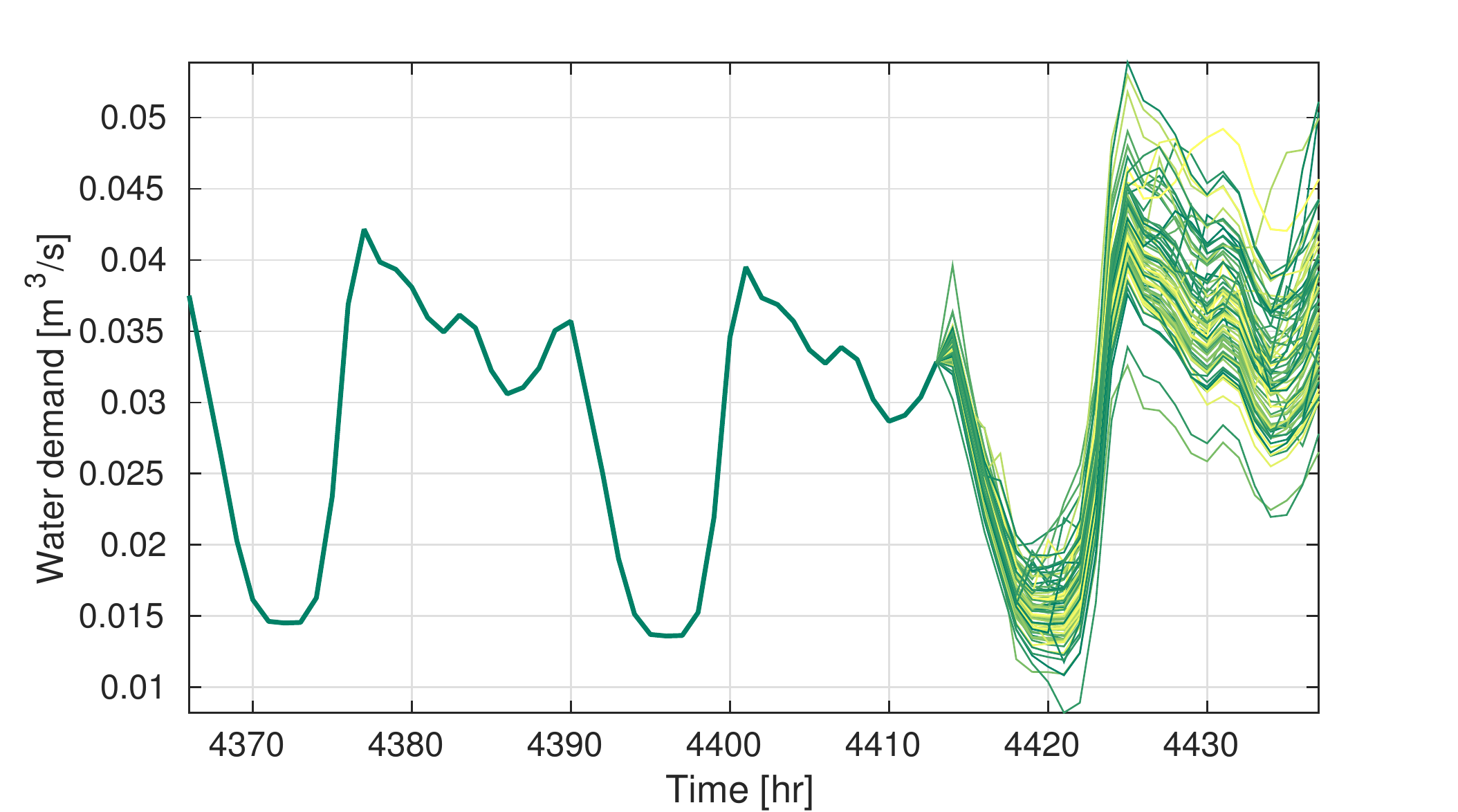}
 \caption{Collection of possible upcoming demands at a given time instant. These results 
 were produced using the SVM model and the data in~\cite{SamGroSop+14}.}
 \label{fig:demadn_realizations}
\end{figure}

\section{Stochastic MPC for DWNs}\label{S.3}
In this section we define the control objectives for the controlled operation of a DWN
and we formulate the stochastic MPC problem.

\subsection{Control objectives}
We define the following three cost functions which reflect our
control objectives. The \textit{economic cost} quantifies the \textit{production} and
 \textit{transportation} cost
 \begin{equation}\label{cost:linear}
  \ell^{w}(u_k,k)=W_{\alpha}(\alpha_1+\alpha_{2,k})'u_k,
 \end{equation}
where the term $\alpha_1' u_k$ is the water production cost, 
$\alpha_{2,k}' u_k$ is the pumping (electricity) cost and 
$W_{\alpha}$ is a positive scaling factor.

The \textit{smooth operation cost} is defined as
\begin{equation}\label{cost:smooth}
 \ell^{\Delta}(\Delta u_k)=\Delta u_k' W_u \Delta u_k,
\end{equation}
where $\Delta u_k=u_k-u_{k-1}$ and $W_u\in\Re^{n_u\times n_u}$ is a symmetric 
positive definite weight matrix. It is introduced to penalize abrupt
switching of the actuators (pumps and valves).

The \textit{safety operation cost} penalizes the drop of water level
in the tanks below a given \textit{safety level}. An elevation
above this safety level ensures that there will be enough water 
in unforeseen cases of unexpectedly high demand
and also maintains a minimum pressure for the flow of water 
in the network. This is given by
\begin{equation}\label{cost:non-smooth}
 \ell^{S}(x_k)=W_x \dist(x_k \mid \mathcal{C}_s),
\end{equation}
where $\dist(x \mid \mathcal{C})=\inf_{y \in \mathcal{C}}\| x-y \|_2$ is the 
distance-to-set function, $\mathcal{C}_s=\{x \mid x \geq x_s \}$, and $x_s\in\Re^{n_x}$ is 
the safety level and $W_x$ is a positive scaling factor.

These cost functions have been used in many MPC formulations in the 
literature~\cite{SamGroSop+14,GroOcam+14,ConPuiCem14}. A comprehensive
discussion on the choice of these cost functions can be found in~\cite{OcFaBa+10}.

The total \textit{stage cost} at a time instant $k$ is the summation of the above 
costs and is given by
\begin{equation}\label{eq:total:cost}
 \ell(x_{k},u_k, u_{k-1},k)=\ell^{w}(u_k,k)+\ell^{\Delta}(\Delta u_k)+\ell^{S}(x_{k}).
\end{equation}


\subsection{SMPC formulation}
We formulate the following stochastic MPC problem with decision 
variables $\pi=\{u_{k+j \mid k},x_{k+j+1 \mid k}\}_{j\in\N_{[0,N-1]}}$
\begin{subequations}\label{eq:SMPC}
\begin{align}
 &V^{\star}(p,q,\hat{\mathbf{d}}_k,k)=\min_{\pi} \expt V(\pi,p,q,k),
 \end{align}
 where $\expt$ is expectation operator and
 \begin{align}\label{eq:total_cost_function}
  V(\pi,p,q,k)=\sum_{j=0}^{N-1}\ell(x_{k+j \mid k},u_{k+j \mid k},u_{k+j-1 \mid k}, k{+}j),
 \end{align}
 subject to the constraints
 \begin{align}
 & x_{k\mid k}=p,\ u_{k-1\mid k}=q,\label{eq:inital_SP}\\
  & x_{k+j+1 \mid k}=Ax_{k+j\mid k}+Bu_{k+j\mid k}+G_dd_{k+j \mid k}(\epsilon_j),\nonumber\\
  &\qquad  j \in \mathbb{N}_{[0,N-1]}, \epsilon_j \in \Omega_j, \label{eq:SMPC:dynamics}\\
  & Eu_{k+j \mid k}+E_dd_{k+j\mid k}(\epsilon_j)=0,  j \in \mathbb{N}_{[0,N-1]},\epsilon_j \in \Omega_j,\label{eq:SMPC:mass-pres}\\
  & x_{\min}\leq x_{k+j\mid k}\leq x_{\max}, j \in \N_{[1,N]},\label{eq:SMPC:state:const}\\
  & u_{\min}\leq u_{k+j\mid k} \leq u_{\max}, j \in \N_{[0,N-1]},
 \end{align}
where we stress out that the decision variables 
$\{u_{k+j \mid k}\}_{j=0}^{j=N-1}$ are required to be causal control 
laws of the form
\begin{align}
u_{k+j \mid k}=\varphi_{k+j\mid k}(p,q, x_{k+j\mid k}, u_{k+j-1\mid k},\boldsymbol{\epsilon}_j).
\end{align}
\end{subequations}

\begin{figure}[ht]
 \centering
 \includegraphics[keepaspectratio=true,width=0.45\textwidth]{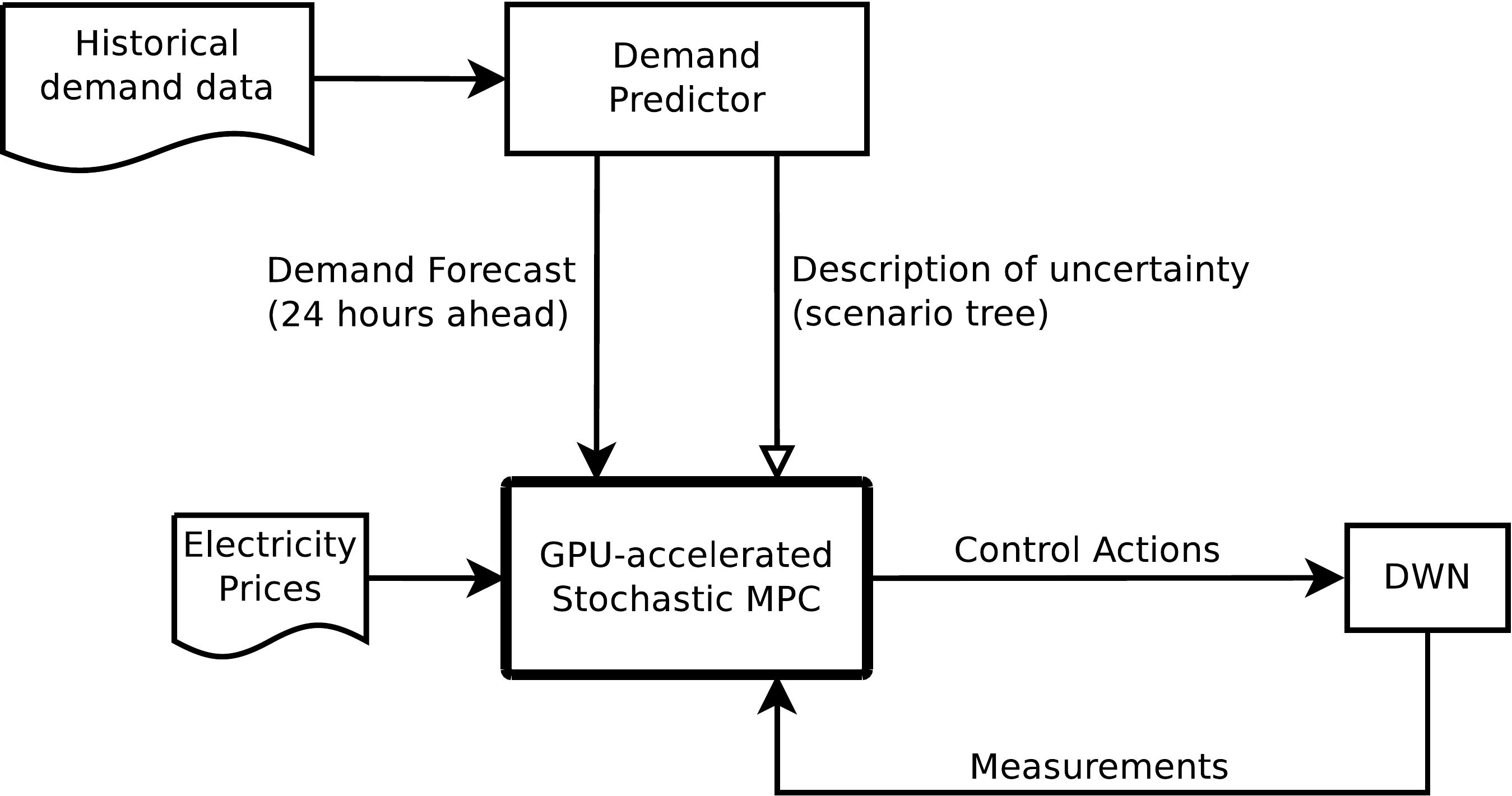}
 \caption{The closed-loop system with the proposed stochastic MPC controller running 
 on a GPU device.}
 \label{fig:smpc_scheme}
\end{figure}

Solving the above problem would involve the evaluation of 
multi-dimensional integrals over an infinite-dimensional space
which is computationally intractable. Hereafter, however, we shall
assume that all $\Omega_j$, for $j\in\N_{[0, N-1]}$, are finite
sets. This assumption will allow us to restate~\eqref{eq:SMPC} as a finite-dimensional
optimization problem.

\subsection{Scenario trees}

\begin{figure}[ht]
 \centering
 \includegraphics[keepaspectratio=true,width=0.48\textwidth]{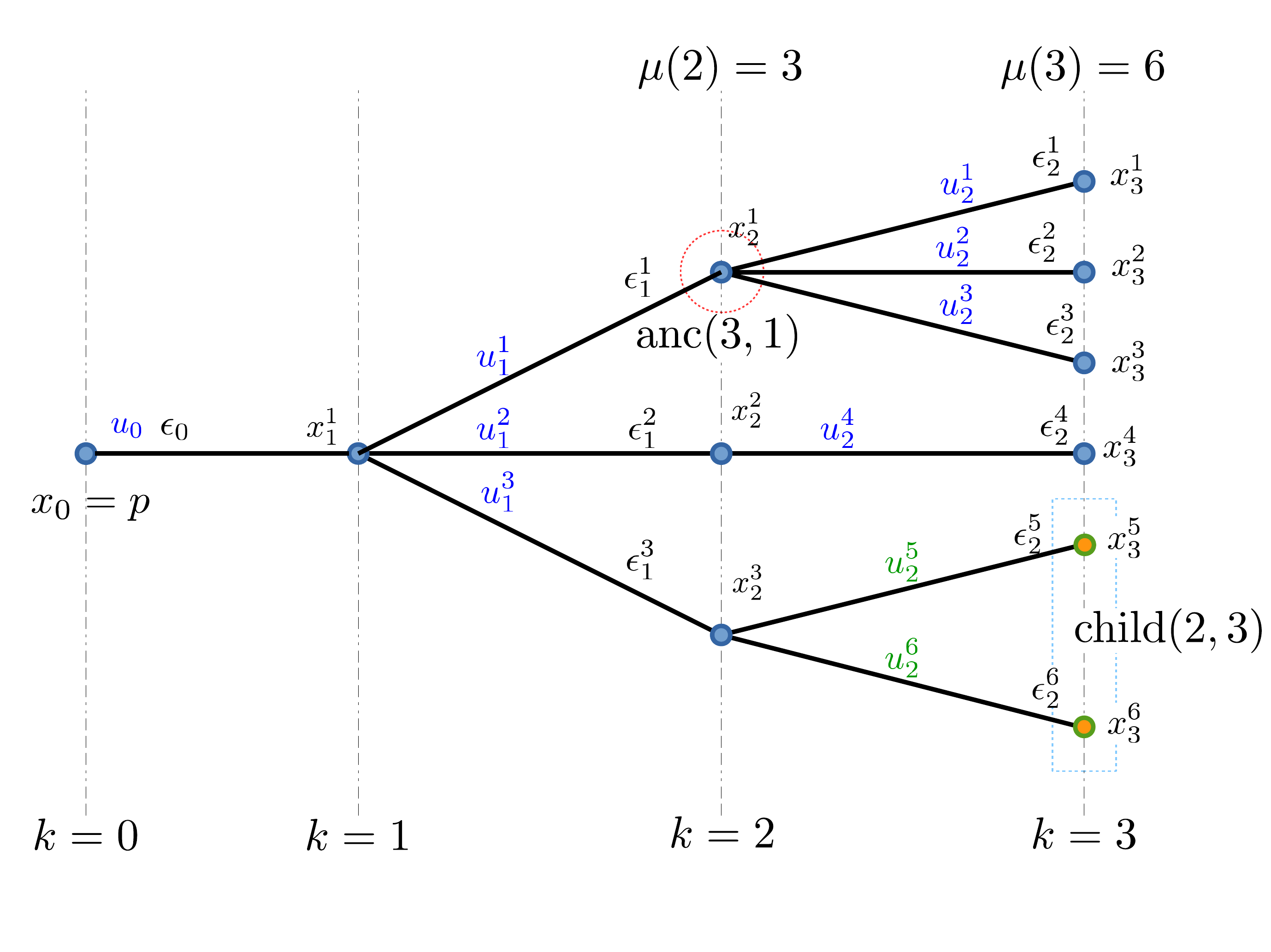}
 \caption{Scenario tree describing the possible evolution of the system state along the 
 prediction horizon: Future control actions are decided in a non-anticipative (causal) fashion; for example  
 $u_1^2$ is decided as a function of $\epsilon_1^2$ but not of any of $\epsilon_2^i$, $i\in\N_{[1,\mu(3)]}$.}
 \label{fig:tree}
\end{figure}

A scenario tree is the structure which naturally follows from the finiteness 
assumption of $\Omega_j$ and is illustrated in Fig.~\ref{fig:tree}. A scenario tree
describes a set of possible future evolutions of the state of the system known as
\textit{scenarios}. Scenario trees can be constructed algorithmically from raw data as 
in~\cite{HeiRom09}. 

The nodes of a scenario tree are partitioned in \textit{stages}.
The (unique) node at stage $k=0$ is called \textit{root} and the nodes
at the last stage are the \textit{leaf nodes} of the tree. 
We denote the number of leaf nodes by $n_s$.
The number of nodes at stage $k$ is denoted by $\mu(k)$ and the total
number of nodes of the tree is denoted by $\mu$.
A path connecting the root node with a leaf node is called 
a \textit{scenario}.
Non-leaf nodes define a set of \textit{children}; at a stage $j\in\N_{[0,N-1]}$
for $i\in\mu(j)$ the set of children of the $i$-th node is denoted by
$\child(j,i)\subseteq \N_{[1,\mu(j+1)]}$. At stage $j\in\N_{[1, N]}$ the
$i$-th node $i\in\N_{[1,\mu(j)]}$ is reachable from a single node
at stage $k-1$ known as its \textit{ancestor} which is denoted by
$\parent(j,i)\in \N_{[1,\mu(j-1)]}$.

The probability of visiting a node $i$ at stage $j$ starting from 
the root is denoted by $p_j^i$. For all for $j \in \N_{N}$ we have that 
$\sum_{i=1}^{\mu(j)}p_j^i=1$  and for all $i\in \N_{[1,\mu(k)]}$ it is 
$ \sum_{l \in \child(j,i)}p_{j+1}^l=p_j^i$.

We define the maximum branching factor at stage $j$, $b_j$, to be
the maximum number of children of the nodes at this stage. The maximum 
branching factor serves as a measure of the complexity of the tree at
a given stage.

\subsection{Reformulation as a finite-dimensional problem} 
We shall now exploit the above tree structure to reformulate the optimal
control problem~\eqref{eq:SMPC} as a finite-dimensional problem.
The water demand, given by~\eqref{eq:demand:model}, is now modeled as
\begin{equation}\label{eq:demand:tree}
 d_{k+j\mid k}^i=\hat{d}_{k+j\mid k}+\epsilon_j^i,
\end{equation}
for all $j \in \N_{[0,N-1]}$ and $i \in \N_{[1,\mu(j+1)]}$.
The input-disturbance coupling~\eqref{eq:SMPC:mass-pres} is then readily 
rewritten as
\begin{equation}\label{eq:scenario:mass-pres}
 Eu_{k+j\mid k}^i+E_d d_{k+j\mid k}^i=0,
\end{equation}
for $j \in \N_{[0,N-1]}$ and $i \in \N_{[1,\mu(j+1)]}$.

The system dynamics is defined across the nodes of the tree by
 \begin{align}
   x_{k+j+1\mid k}^l  &=  A x_{k+j\mid k}^i  +  B u_{k+j\mid k}^l  +  G_d d_{k+j\mid k}^l,
 \end{align}
for $j \in \N_{[0,N-1]}$, $i \in \N_{[1,\mu(j)]}$ and $l \in \child(j,i)$, or, alternatively,
 \begin{align}\label{eq:scenario:dynamics}
   x_{k+j+1\mid k}^i  &=  A x_{k+j\mid k}^{\parent(j+1,i)}  +  B u_{k+j\mid k}^i  +  G_d d_{k+j\mid k}^i,
 \end{align}   
for $j \in \N_{[0,N-1]}$ and $i\in\N_{[1,\mu(j+1)]}$.
   
Now the expectation of the objective function~\eqref{eq:total_cost_function} can be derived
as a summation across the tree nodes
\begin{align}
 &\expt V(\pi,p,q,k)=\notag\\
 &\quad\sum_{j=0}^{N-1}\sum_{i=1}^{\mu(j)}p_j^i\ell(x_{k+j\mid k}^i, u_{k+j\mid k}^i, u_{k+j-1 \mid k}^{\parent(j,i)},k+j), 
\end{align}
where $x_{k\mid k}^1=p$ and $u_{k-1\mid k}=q$.

In order to guarantee the recursive feasibility of the control problem, 
the state constraints~\eqref{eq:SMPC:state:const} are converted into 
\textit{soft constraints}, that is, they are replaced by a penalty of the 
form
\begin{equation}\label{cost:distance}
 \ell^{d}(x)=\gamma_d\dist(x,\mathcal{C}_1),
\end{equation}
where $\gamma_d$ is a positive penalty factor and $\mathcal{C}_1=\{x \mid 
x_{\min} \leq x \leq x_{\max}\}$. 
Using this penalty, we construct the 
\textit{soft state constraint penalty}
\begin{align}
 V_{s}(\pi,p) = \sum_{j=0}^{N}\sum_{i=1}^{\mu(j)}\ell^d(x_{k+j\mid k}^{i}).
\end{align}

The modified, soft-constrained, SMPC problem can be now written as
\begin{subequations}\label{eq:Sc-MPC}
\begin{align}
 \tilde{V}^{\star}(p,q,\hat{\mathbf{d}},k)=&\min_{\pi} \expt V(\pi,p,q,k) +V_{s}(\pi,p),
 \end{align}
subject to 
 \begin{align}
 & x_{k\mid k}^1=p,\ u_{k-1\mid k}=q,\\
  & u_{\min}\leq u_{k+j\mid k}^i \leq u_{\max}, j \in \N_{[0,N-1]}, i \in \N_{[1,\mu(j)]},\label{const:control}
 \end{align}
\end{subequations}
and system equations~\eqref{eq:scenario:mass-pres} and~\eqref{eq:scenario:dynamics}.

\section{Solution of the stochastic optimal control problem}\label{S.5}

In this section we extend the GPU-based proximal gradient method proposed
in~\cite{SaPaBe+15} to solve the SMPC problem~\eqref{eq:Sc-MPC}.
For ease of notation we will focus on the solution of the SMPC problem at
$k=0$ and denote $x_{j\mid 0}=x_j$, $u_{j\mid 0}=u_j$, $\hat{d}_{j\mid 0}=\hat{d}_{j}$.

\subsection{Proximal gradient algorithm}
For a closed, proper extended-real valued function 
$g:\Re^n \rightarrow \bar{\Re}$, we define its 
\textit{proximal operator} with parameter $\gamma>0$, 
$\prox_{\gamma g}: \Re^n \rightarrow \Re^n$ as~\cite{ParBoy13}
\begin{equation}
 \prox_{\gamma g}(v)= \argmin_{x\in \Re^n}\bigg\{ g(x)+\frac{1}{2\gamma}\| x-v\|_2^2\bigg\}.
\end{equation}
The proximal operator of many functions is available in closed form~\cite{ComPes10,ParBoy13}.
When $g$ is given in a \textit{separable sum} form, that is
\begin{subequations}\label{eq:separable-sum-property}
\begin{align}
 g(x) = \sum_{i=1}^{\kappa}g_i(x_i),
\end{align}
then,
\begin{align}
 (\prox_{\gamma g}(v))_i = \prox_{\gamma g_i}(v_i).
\end{align}
\end{subequations}
This is known as the \textit{separable sum property} of the proximal operator.

Let $z\in\Re^{n_z}$ be a vector encompassing all states $x_j^i$ 
for $j\in\N_{[0, N]}$ and $i\in \N_{[1,\mu(j)]}$ and inputs $u_j^i$ 
for $j\in\N_{[0, N-1]}$, $i\in \N_{[1,\mu(j+1)]}$; this is the 
decision variable of problem~\eqref{eq:Sc-MPC}.

Let $f:\Re^{n_z}\to \bar{\Re}$ be defined as 
\begin{align}
 f(z)  &= \sum_{j=0}^{N-1}\sum_{i=1}^{\mu(j)}  p_j^i(\ell^w(u_j^i)+\ell^{\Delta}(\Delta u_j^i))\nonumber\\
                      & +\delta(u_j^i|\Phi_1(d_j^i)) \notag\\
                      &+ \delta(x_{j+1}^{i}, u_j^i, x_j^{\parent(j+1,i)}|\Phi_2(d_j^i)),\label{eq:fz}
\end{align}
where $\Delta u_j^i= u_j^i- u_{j-1}^{\parent(j,i)}$ and 
$\Phi_1(d)$ is the affine subspace of $\Re^{n_u}$ induced by~\eqref{eq:scenario:mass-pres},
that is
\begin{equation}\label{eq:Phi_1}
\Phi_1(d) = \{u:Eu + E_d d = 0\},
\end{equation}
and $\Phi_2(d)$ is the affine subspace of $\Re^{2n_x+n_u}$ defined by the system
dynamics~\eqref{eq:scenario:dynamics}
\begin{equation}\label{eq:Phi_2}
\Phi_2(d)=\{(x_{k+1}, x_k, u): x_{k+1}=Ax_k + Bu + G_d d\}.
\end{equation}

We define the auxiliary variables $\varsigma$ and $\zeta$ which stand 
for \textit{copies} of the state variables $x_j^i$  --- that is 
$\varsigma_j^i = \zeta_j^i = x_j^i$ ---
and the auxiliary variable $\psi$ 
which is a copy of input variables $\psi_j^i = u_j^i$. The reason for the introduction 
of these variables will be clarified in Section~\ref{sec:dual-iterate}.

We introduce the variable $t=(\varsigma,\zeta,\psi)\in\Re^{n_t}$  
and define an extended real valued function $g:\Re^{n_t}\to\bar{\Re}$ as
\begin{align}
g(t) &= \sum_{j=0}^{N-1} \sum_{i=1}^{\mu(j)}\ell^{S}(\varsigma_{j+1}^i)
+\ell^{d}(\zeta_{j+1}^i)+ \delta(\psi_j^i| \mathcal{U}),\label{eq:gzeta}
\end{align}
where
$\mathcal{U}=\{\psi\in\Re^{n_u}: u_{\min}\leq \psi \leq u_{\max}\}$.

Now the finite-dimensional optimization problem~\eqref{eq:Sc-MPC} can be written as:
\begin{subequations}\label{eq:sc_prox}
\begin{align}
 \tilde{V}^{\star}  &=  \min_{z,t} f(z)+g(t), \\
 \text{s.t. } &Hz=t,
\end{align}
\end{subequations}
where 
\begin{equation}\label{optim:H:operator}
 H=\left[ \begin{array}{cc}
          I_{n_x} & 0\\
          I_{n_x} & 0\\
          0 & I_{n_u}
         \end{array} \right].
\end{equation}

The Fenchel dual of~\eqref{eq:sc_prox} is written as~\cite[Corol.~{31.2.1}]{Roc72}:
\begin{align}\label{optim:prob}
 \tilde{D}^{\star}=\min_{y} f^{*}(-H' y)+g^{*}(y),
\end{align}
where $y$ is the dual variable. 
The dual variable $y$ can be partitioned as $y=(\tilde\varsigma_{j}^{i}, \tilde\zeta_j^i, \tilde\psi_j^i)$,
where $\tilde\varsigma_{j}^{i}$, $\tilde\zeta_j^i$ and $\tilde\psi_j^i$ are the dual variables corresponding
to $\varsigma_{j}^{i}$, $\zeta_j^i$ and $\psi_j^i$ respectively. We also 
define the auxiliary variable of state copies $\tilde{\xi}_j^i \dfn (\tilde{\varsigma}_j^i,\tilde{\zeta}_j^i)$.

According to~\cite[Thm.~{11.42}]{RocWets09}, since function $f(z)+g(Hz)$ is proper, convex and piecewise linear-quadratic, then
the primal problem~\eqref{eq:sc_prox} is feasible whenever the dual problem~\eqref{optim:prob} is feasible
and, furthermore, strong duality holds, i.e., $\tilde{V}^{\star}=\tilde{D}^{\star}$. Moreover, the optimal solution of~\eqref{eq:sc_prox}
is given by $z^\star = \nabla f^{\ast}(-H'y^\star)$ where $y^\star$ is \textit{any} solution of~\eqref{optim:prob}. 
Applying~\cite[Prop.~{12.60}]{RocWets09} to $f^*$ and since $f$ is lower semi-continuous, proper and 
$\sigma$-strongly convex --- as shown at the end of Appendix~\ref{App:cost update} --- 
its conjugate $f^{\ast}$ has Lipschitz-continuous gradient 
with a constant $1/\sigma$.

An accelerated version of proximal-gradient method which was first proposed 
by Nesterov in~\cite{Nes83} is applied to the dual problem. This leads to the 
following algorithm
\begin{subequations}\label{eq:apg-iteration}
 \begin{align}
  w^\nu &=y^\nu+\theta_\nu(\theta_{\nu-1}^{-1}-1)(y^\nu-y^{\nu-1}),\label{optim:acceleration}\\
  z^\nu&=\argmin_{z}\{\langle z,H^{\prime} w^\nu\rangle +f(z)\},\label{optim:primal:update}\\
  t^\nu&=\prox_{\lambda^{-1}g}(\lambda^{-1}w^{\nu}+Hz^\nu),\label{optim:dual:update1}\\
  y^{\nu+1}&=w^\nu+\lambda(Hz^v-t^v),\label{optim:dual:update2}\\
  \theta_{\nu+1}&=\frac{1}{2}\bigg(\sqrt{\theta_\nu^4+4\theta_\nu^2}-\theta_\nu^2\bigg)\label{optim:theta:update},
 \end{align}
\end{subequations}
starting from a dual-feasible vector $y^0=y^{-1}=0$ and $\theta_0=\theta_{-1}=1$.

In the first step~\eqref{optim:acceleration} we compute an extrapolation of the dual vector. 
In the second step~\eqref{optim:primal:update} we calculate the dual gradient, 
that is $z^{\nu}=\nabla f^{\ast}(-H' w^{\nu})$, 
at the extrapolated dual vector using the conjugate subgradient theorem~\cite[Thm.~{23.5}]{Roc72}.
The third step comprises of~\eqref{optim:dual:update1},~\eqref{optim:dual:update2} where we
update the dual vector $y$ and in the final step of the algorithm we compute the scalar $\theta_{\nu}$
which is used in the extrapolation step.

This algorithm has a convergence rate of $\mathcal{O}(1/\nu^2)$ for the 
dual iterates as well as for the ergodic primal iterate defined through 
the recursion $\bar{z}^\nu=(1-\theta_{\nu})\bar{z}^{(\nu-1)}+\theta_{\nu}z^{\nu}$, i.e.,
a weighted average of the primal iterates~\cite{PanBem14}.

\subsection{Computation of primal iterate}
The most critical step in the algorithm is the computation of $z^\nu$
which accounts for most of the computation time required by each iteration. This step
boils down to the solution of an unconstrained optimization problem by means of
dynamic programming where certain matrices (which are independent of $w^\nu$) can be computed
once before we run the algorithm to facilitate the online computations.
These are (i) the vectors $ \beta_j^i, \hat{u}_j^i,e_j^i$ which are associated with the update
of the time-varying cost (see Appendix~\ref{App:cost update}) and 
(ii) the matrices $\Lambda, \Phi, \Psi,\bar{B}$ (see Appendix~\ref{App:Factor step}).
The latter are referred to as the \textit{factor step} of the algorithm and matrices
$\Lambda, \Phi, \Psi$ and $\bar{B}$ are independent of the complexity of the scenario tree. 

The computation of $z^{\nu}$ at each iteration of the algorithm requires the computation
of the aforementioned matrices and is computed using Algorithm~\ref{algo:solve_dp} to which 
we refer as the \textit{solve step}. Computations involved in the solve step are merely
matrix-vector multiplications. As the algorithm traverses the nodes of the 
scenario tree stage-wise backwards (from stage $N-1$ to stage $0$), 
computations across the nodes at a given stage can be 
performed in parallel. Hardware such as GPUs which enable us to parallelizable
such operations lead to a great speed-up as we demonstrate in Section~\ref{S.6}.

\begin{algorithm}[ht]
\caption{Solve step}
\label{algo:solve_dp}
\begin{algorithmic}
\REQUIRE Output of the factor step (See Appendices~\ref{App:cost update} and~\ref{App:Factor step}), i.e., 
$\Lambda, \Phi, \Psi, \bar{B},\hat{u}_j^i,\beta_{j}^{i}, e_{j}^{i}$, $p$, $q$ and 
$w^\nu = (\tilde{\varsigma}_j^i, \tilde{\zeta}_j^i, \tilde{\psi}_j^i)$.
\STATE $q_{N}^i \set 0$, and $r_{N}^i\set0,\forall i\in\N_{[1,n_s]}$,
\FOR{$j=N-1,\ldots,0$}
   \FOR[in parallel]{$i=1,\ldots,\nodes(k)$}
      \STATE $\sigma_j^l\set r_{j+1}^l+\beta_j^l, \forall l \in \child(j,i)$
      \STATE $v_j^l\set \frac{1}{2p_j^l}\left(\Phi_j^l(\tilde{\xi}_j^l+q_{j+1}^l)  +   \Psi_j^l\tilde{\psi}_j^l+\Lambda_j^l\sigma_j^l
                    \right),$ \\ \quad $\forall l \in \child(j,i)$
      \STATE $r_j^i\set \sum_{l \in \child(j,i)} \sigma_j^l  +  \bar{B}^{\prime}(\tilde{\xi}_j^l+q_{j+1}^l)   +    L\tilde{\psi}_j^l$
      \STATE $q_j^i\set A^{\prime}\sum_{l \in \child(j,i)}\tilde{\xi}_j^l+q_{j+1}^l$
   \ENDFOR 
\ENDFOR 
\STATE $x_0^1\set p$, $u_{-1}\set q$, 
\FOR{$j=0,\ldots,N-1$}
   \FOR[in parallel]{$i=1,\ldots,\nodes(k)$}
      \STATE $v_j^i\set v_{j-1}^{\parent(j,i)}+v_j^i$ 
      \STATE $u_j^i\set Lv_j^i+\hat{u}_j^i$ 
      \STATE $x_{j+1}^i\set Ax_j^{\parent(j,i)}+\bar{B}v_j^i+e_{j}^i$
   \ENDFOR
\ENDFOR 
\RETURN $\{x_j^i\}_{j=1}^{N}$, $\{u_j^i\}_{j=0}^{N-1}$
\end{algorithmic}
\end{algorithm}

\subsection{Computation of dual iterate}\label{sec:dual-iterate}
Function $g$ given in~\eqref{eq:gzeta} is given in the form of a separable sum
\begin{equation}\label{optim:g}
g(t) = g(\varsigma,\zeta, \psi) = g_1(\varsigma)+g_2(\zeta)+g_3(\psi),
\end{equation}
where 
\begin{subequations}
 \begin{align}
 g_1(\varsigma)&= \sum_{j=0}^{N-1}\sum_{i=1}^{\mu(j)} \ell^S(\varsigma_{j+1}^i),\label{optim:prox:safe}\\
 g_2(\zeta)    &= \sum_{j=0}^{N-1}\sum_{i=1}^{\mu(j)} \ell^d(\zeta_{j+1}^i),\label{optim:prox:x}\\
 g_3(\psi)     &= \sum_{j=0}^{N-1}\sum_{i=1}^{\mu(j)} \delta(\psi_j^i\mid \mathcal{U})\label{optim:prox:u}.
\end{align}
\end{subequations}

Functions $g_1(\cdot)$ and $g_2(\cdot)$ are in turn separable sums of 
distance functions from a set and $g_3(\cdot)$ is an indicator 
function. Their proximal mappings can be easily computed as in Appendix~\ref{App:Proximal}
and essentially are element-wise operations on the vector $t$ that can be fully parallelized.

\subsection{Preconditioning and choice of $\lambda$}\label{sec:preconditioning}

First-order methods are known to be sensitive to scaling and preconditioning can remarkably improve their 
convergence rate. Various preconditioning method such as~\cite{GisBoy15,Bradley10} have 
been proposed in the literature. Here, we employ a simple diagonal preconditioning which 
consists in computing a diagonal matrix $\tilde{H}_D$ with positive diagonal entries which approximates 
the dual Hessian $H_D$ and use $\tilde{H}^{-1/2}_D$ to scale the dual vector~\cite[2.3.1]{Ber99}. 
Since the uncertainty does not affect the dual Hessian, we take this preconditioning matrix for 
a single branch of the scenario tree and use it to scale all dual variables.

In a similar way, we compute the parameter $\lambda$. We choose $\lambda = 1/L_{H_D}$ where $L_{H_D}$ is
the Lipschitz constant of the dual gradient which is computed as $\|H\|^2/\sigma$ as in~\cite{Ber99}. 
It again suffices to perform the computation for a single branch of the scenario tree.

\subsection{Termination}
The termination conditions for the above algorithm are based on the ones
provided in~\cite{PanBem14}. However, rather than checking these conditions
at every iteration, we perform always a fixed number of iterations which 
is dictated by the sampling time. We may then check the quality of the solution
\textit{a posteriori} in terms of the duality gap and the term $\|Hz^\nu-t^\nu\|_{\infty}$.

\section{Case study: The Barcelona DWN}\label{S.6}
\begin{figure}
\centering
 \includegraphics[keepaspectratio=true,width=0.45\textwidth]{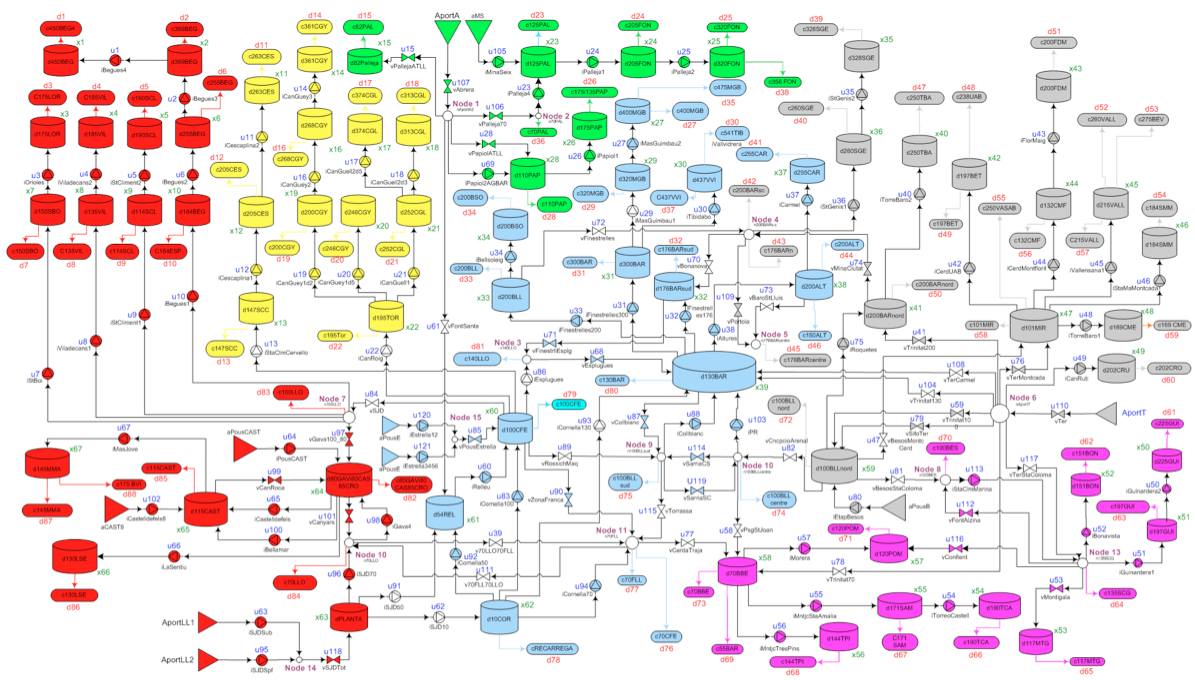}
 \caption{Structure of the DWN of Barcelona.}\label{fig:dwn}
\end{figure}

We now apply the proposed control methodology to the drinking water network of
the city of Barcelona using the data found in~\cite{GroOcam+14,SamGroSop+14}.
The topology of the network is presented in Figure~\ref{fig:dwn}.
The system model consists of $63$ states corresponding to the level of water in each 
tank, $114$ control inputs which are pumping actions and valve positions, $88$ demand nodes
and $17$ junctions. The prediction horizon is $N=24$ with sampling time of $1$ hour. 
The future demands are predicted using the SVM time series model developed in~\cite{SamGroSop+14}. 

\subsection{Performance of GPU-accelerated algorithm}
Accelerated proximal gradient (APG) was implemented in CUDA-C v6.0 and the matrix-vector computations were performed using 
cuBLAS. We compared the GPU-based implementation with the interior-point solver of Gurobi.
Active-set algorithms exhibited very poor performance and we did not include the respective 
results.

All computations on CPU were performed on a $4\times2.60\mathrm{GHz}$ Intel i5 machine 
with $8\mathrm{GB}$ of RAM running 64-bit Ubuntu 
v14.04 and GPU-based computations were carried out on a NVIDIA Tesla C2075. 

\begin{figure}[ht]
 \centering
 \includegraphics[keepaspectratio=true,width=0.48\textwidth]{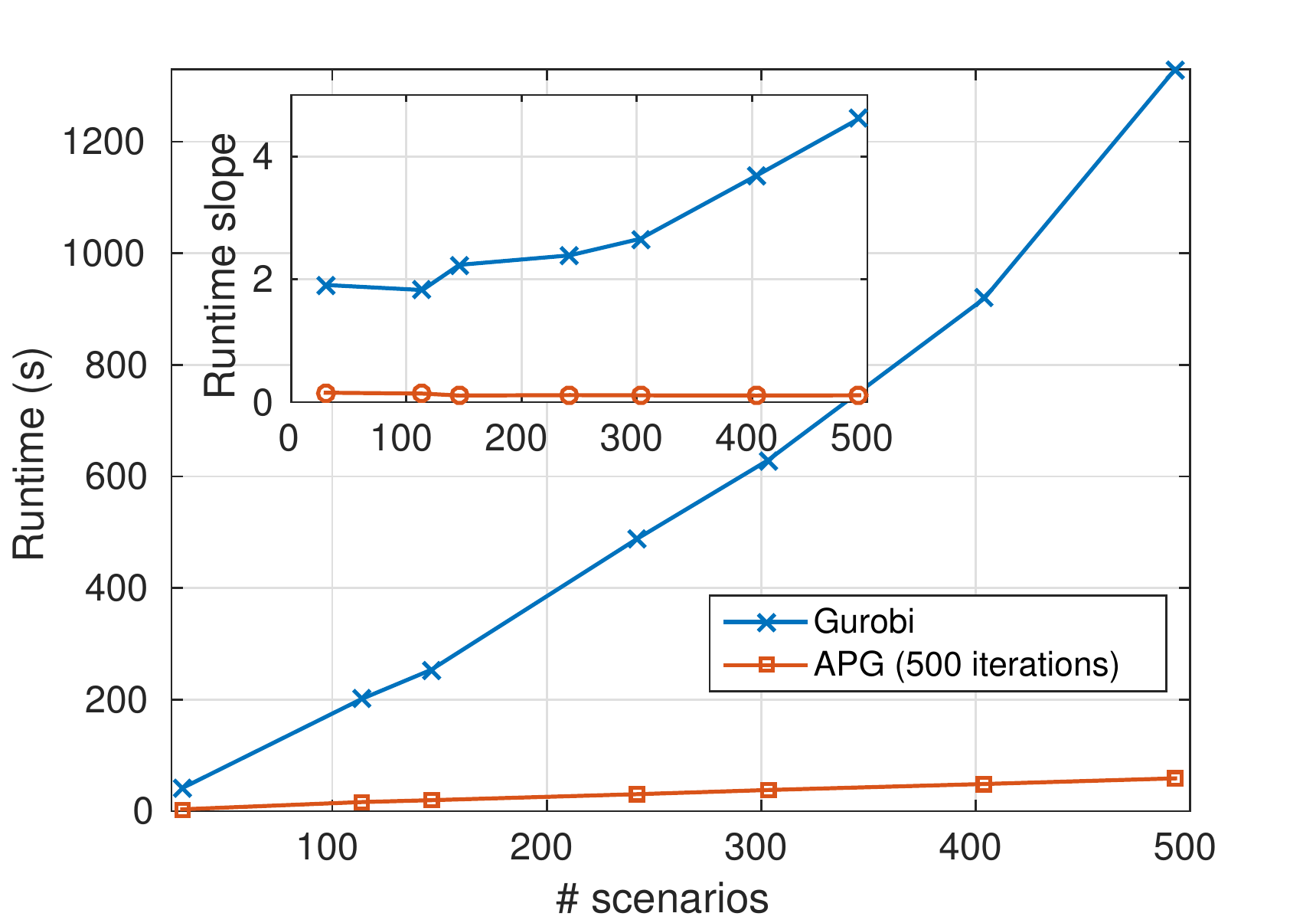}
 \caption{Runtime of the CUDA implementation against the number of scenarios considered 
 in the optimization problem. Comparison with the runtimes of Gurobi.}
 \label{fig:runtimes}
\end{figure}

The dependence of the computational time on the size of the scenario tree 
is reported in the Figure~\ref{fig:runtimes} where it can be noticed that 
there is speed-up of $10\times$ to $25\times$ in the computational times with 
CUDA-APG compared to Gurobi. Furthermore, the speed-up increases with the 
number of scenarios. 

The optimization problems we are solving here are of noticeably large
size. Indicatively, the scenario tree with $493$ scenarios counts 
approximately $2.52$ million dual decision variables ($1.86$ million primal variables)
and while Gurobi requires $1329s$ to solve it, our CUDA implementation solves it in $58.8s$; 
this corresponds to a  speed-up of $22.6\times$.

%
%
%
In all of our simulations we obtained a sequence of control actions across the tree nodes $U^\star_{apg}=\{u_j^i\}$ which 
was, element-wise, within $\pm 0.029m^3/s$ ($1.9\%$) from the solution produced by Gurobi.
The maximum primal residual was $\|Hx-z\|_{\infty}=1.7446$. Moreover, we should note that the control action $u_0^\star$ computed
by APG with $500$ iterations was consistently within $\pm 0.0025 m^3/s$ ($0.08\%$) from the Gurobi solution.
Given that only $u_0^\star$ is applied to the system while all other control actions $u_j^i$ for $j\in\N_{[1,N-1]}$ 
and $i\in\N_{[1, \mu(j)]}$ are discarded, $500$ iterations are well sufficient for convergence.

\subsection{Closed-loop performance}
In this section we analyse the performance of SMPC with different scenario-trees. This analysis is 
carried for a period of 7 days ($H_s=168$) from 1$^\text{st}$ to 8$^\text{th}$ July 2007.
Here, we compare the operational cost and the quality of service of various
scenario-tree structures. 

The weighting matrices in the operational cost are chosen as 
$W_{\alpha}=2\cdot 10^4$, $W_{u}=10^5\cdot I$ and $W_x=10^7$, respectively
and $\gamma_d = 5\cdot 10^{7}$.
The demand is predicted using SVM model presented in~\cite{SamGroSop+14}. 
The steps involved in SMPC using GPU based APG in closed-loop is 
summarized in Algorithm~\ref{algo:closed:loop}.

\begin{algorithm}[ht]
\caption{Closed-loop of DWN with SMPC with proximal-operator}
\label{algo:closed:loop}
\begin{algorithmic}
\REQUIRE Scenario-tree, current state measurement $x_0$ and previous control $u_{-1}$.
\STATE Compute $\Lambda$, $\Phi$, $\Psi$ and $\bar{B}$ as in Appendix~\ref{App:Factor step}
\STATE Precondition the original optimization problem and compute $\lambda$ as in Section~\ref{sec:preconditioning}.       
\LOOP
  \STATE \textit{Step 1.} Predict the future water demands $\hat{\mathbf{d}}_{k}$ using current and past demand data.
  \STATE \textit{Step 2.} Compute $\hat{u}_j^i$, $\beta_j^i,e_j^i$ as in Appendix~\ref{App:cost update}.
  \STATE \textit{Step 3.} Solve the optimization problem using APG on GPU using iteration~\eqref{eq:apg-iteration} and
	Algorithm~\ref{algo:solve_dp}.
  \STATE \textit{Step 4.} Apply $u_0^1$ to the system, update the previous $u_{-1}=u_0^1$
\ENDLOOP
\end{algorithmic}
\end{algorithm}

For the performance assessment of the proposed control methodology we used various 
controllers summarized in Table~\ref{table:controllers}. 
The corresponding computational times are presented in Figure~\ref{fig:runtimes}.  
\begin{table}[ht]
\centering
 \begin{tabular}{c|c|c|c|c}
 Controller & $b_k$ & scenarios & primal variables & dual variables \\
 \hline
 CE-MPC     & $1$   & $1$   & $4248$     & $5760$\\
 \hline
 SMPC$_1$ & $[3,2]$   & $6$   & $24072$  & $32540$\\
 \hline
 SMPC$_2$ & $[6,5]$   & $30$  & $118059$ & $160080$\\
 \hline
 SMPC$_3$ & $[6,5,5]$ & $114$ & $430287$ & $583440$\\
 \hline
 SMPC$_4$ & $[8,5,5]$ & $146$ & $551355$ & $747600$ \\
 \hline
 SMPC$_5$ & $[10,8,5]$& $242$ & $915621$ & $1241520$\\
\hline
SMPC$_6$  & $[12,8,5]$& $303$ & $1145544$ & $1553280$\\
\hline
SMPC$_7$  & $[12,8,8]$& $404$& $1520961$ & $2062320$\\
\hline
SMPC$_8$  & $[12,10,8]$& $493$& $1856022$ & $2516640$\\
\hline
 \end{tabular}
 \caption{Various controllers used to assess the closed-loop performance of the proposed methodology.
 The numbers in the bracket denote the first \textit{maximum branching factors,} $b_j$, of the scenario tree while all
 subsequent branching factors are assumed to be equal to $1$.}\label{table:controllers}
\end{table} 

To assess the performance of closed-loop operation of the SMPC-controlled
network we used the \textit{key performance indicators} (KPIs)
reported in~\cite{AlBa+06,GroMaeOca14}. For a simulation time length $H_s$
the performance indicators are computed by
\begin{subequations}
 \begin{align}
 &\mathrm{KPI}_E=\frac{1}{H_s}\sum_{k=1}^{H_s}(\alpha_1 + \alpha_{2,k})'|u_k|,\\
 &\mathrm{KPI}_{\Delta U}=\frac{1}{H_s}\sum_{k=1}^{H_s}\|\Delta u_{k}\|^2,\\
 &\mathrm{KPI}_{S}=\sum_{k=1}^{H_s}\|[x_{s}-x_{k}]_+\|_1,\\
 &\mathrm{KPI}_R=\frac{ \|x_{s}\|_1}
                  {\frac{1}{H_s}\sum_{k=1}^{H_s}  \|x_{k}\|_1 }\times 100 \%.
\end{align}
\end{subequations}
$\mathrm{KPI}_E$ is the average economic cost, 
$\mathrm{KPI}_{\Delta U}$ measures the average smoothness of the control actions, 
$\mathrm{KPI}_{S}$ corresponds to the total amount of water used from 
storage and $\mathrm{KPI}_R$ is the percentage of the safety volume $x_s$ 
contained into the average volume of water.

\begin{table}[ht]
\centering
 \begin{tabular}{c|c|c|c|c}
 
 Controller & $\mathrm{KPI}_E$ & $\mathrm{KPI}_{\Delta U}$ & $\mathrm{KPI}_{S}$ & $\mathrm{KPI}_{R}$\\
 \hline
 CE-MPC   &  $\mathbf{1801.4}$  & $\mathbf{0.2737}$   	& $\mathbf{6507.7}$   & $64.89\%$\\
 \hline
 SMPC$_1$   & $1633.5$            & $0.3896$            & $1753.7$           & $\mathbf{67.96\%}$\\
 \hline
 SMPC$_2$   & $\mathbf{1549.7}$   & $0.4652$   		& $2264.0$           & $61.81\%$\\
 \hline
 SMPC$_3$   & $1574.0$            & $0.4135$            & $1360.0$           & $49.65\%$\\
 \hline
 SMPC$_4$   & $1583.2$            & $0.4088$            & $885.7 $           & $48.13\%$\\
 \hline
 SMPC$_5$   & $1597.3$            & $0.4470$            & $508.5$   	     & $46.05\%$ \\
 \hline
 SMPC$_6$   & $1606.3$            & $\mathbf{0.4878}$            & $\mathbf{302.3}$   & $\mathbf{44.93\%}$ \\
 \end{tabular}
  \caption{KPIs for performance analysis of the DWN with different controllers. The lowest and the highest in each of 
  the indicator is highlighted. The economical benefit and risk is presented with terms of number of scenarios  
  }\label{table:KPI}
\end{table}

\subsubsection{Risk vs Economic utility}
\begin{figure}[ht]
 \centering
 \includegraphics[keepaspectratio=true,width=0.48\textwidth]{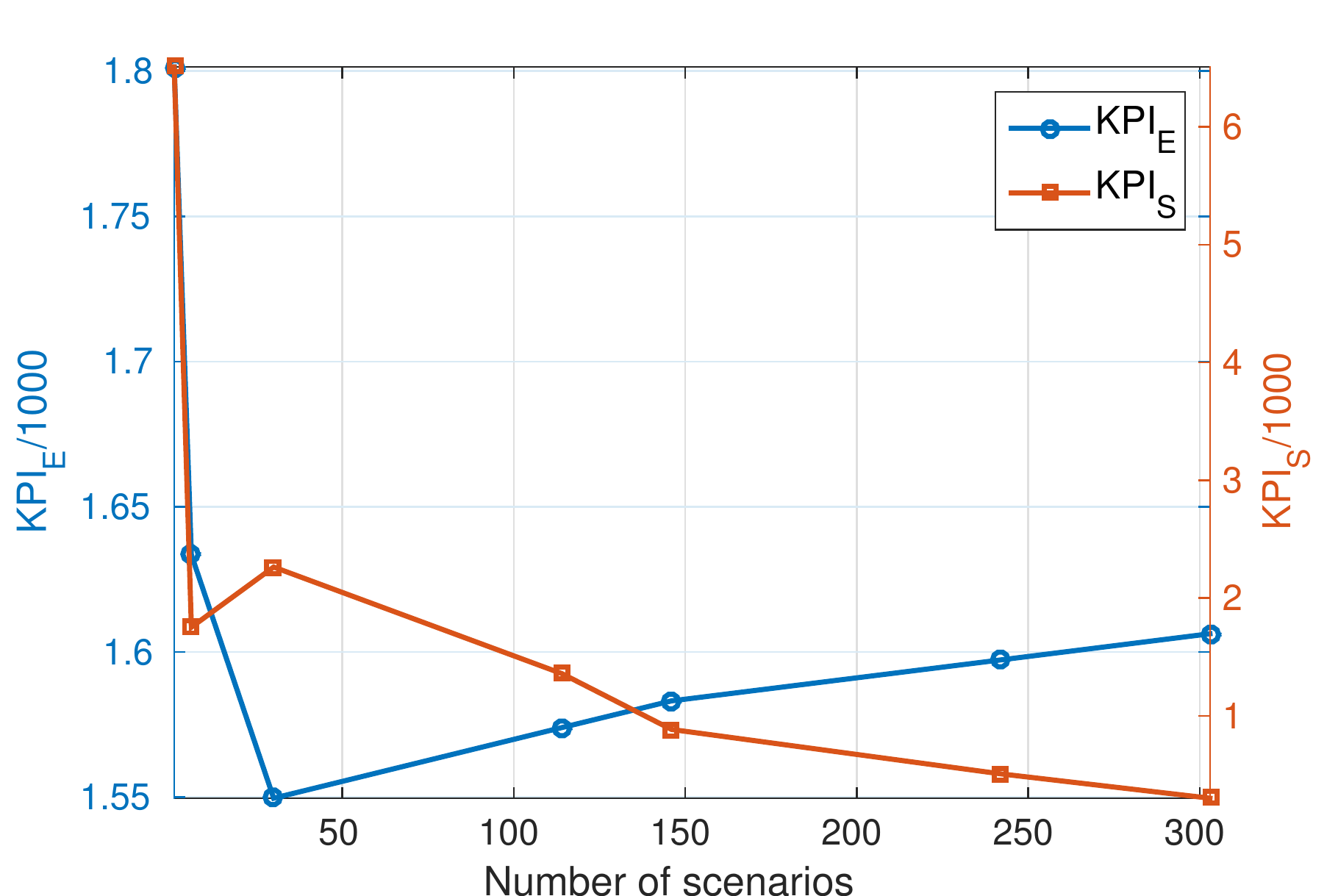}
 \caption{The figure shows the trade-off between risk and economic utility 
 in terms of scenarios. The KPI$_E$ represent the economical 
 utility and KPI$_S$ shows the risk of violation.}
 \label{fig:kpis}
\end{figure}

Figure~\ref{fig:kpis} illustrates the trade-off between economic and safe operation:
The more scenarios we use to describe the distribution of demand prediction 
error, the safer the closed-loop operation becomes as it is reflected by the 
decrease of KPI$_S$. Stochastic MPC leads to a significant decrease of economic
cost compared to the certainty-equivalence approach, however, the safer we 
require the operation to be, the higher the operating cost we should expect.

\subsubsection{Quality of service}
A measure of the reliability and quality-of-service of the 
network is $\mathrm{KPI}_S$ which reflects the tendency of 
water levels to drop under the safety storage levels.
As expected, the CE-MPC controller leads to the most unsafe 
operation, whereas SMPC$_6$ leads to the lowest value.

\subsubsection{Network utility}
Network utility is defined as the ability to utilize the water in 
the tanks to meet the demands rather than pumping additional water
and is quantified by $\mathrm{KPI}_R$. In Table~\ref{table:KPI}, we see the dependence 
of $\mathrm{KPI}_R$ on the number of scenarios of the tree. 
$\mathrm{KPI}_R$ remains always within reasonable limits; on average we 
operate away from the safety storage limit. The decrease in $\mathrm{KPI}_R$
on may observe is because as more scenarios are employed, the more accurate the 
representation of uncertainty becomes and the system does not need to operate,
on average, too far away from $x_s$.

\subsubsection{Smooth operation}
We may notice that the introduction of more scenarios results in 
an increase in $\mathrm{KPI}_{\Delta U}$. Then, the controller becomes
more responsive to accommodate the need for a less risky operation, 
although the value of $\mathrm{KPI}_{\Delta U}$ is not greatly affected
by number of scenarios.

\begin{figure}[ht]
 \centering
 \includegraphics[width=0.49\textwidth,keepaspectratio=true]{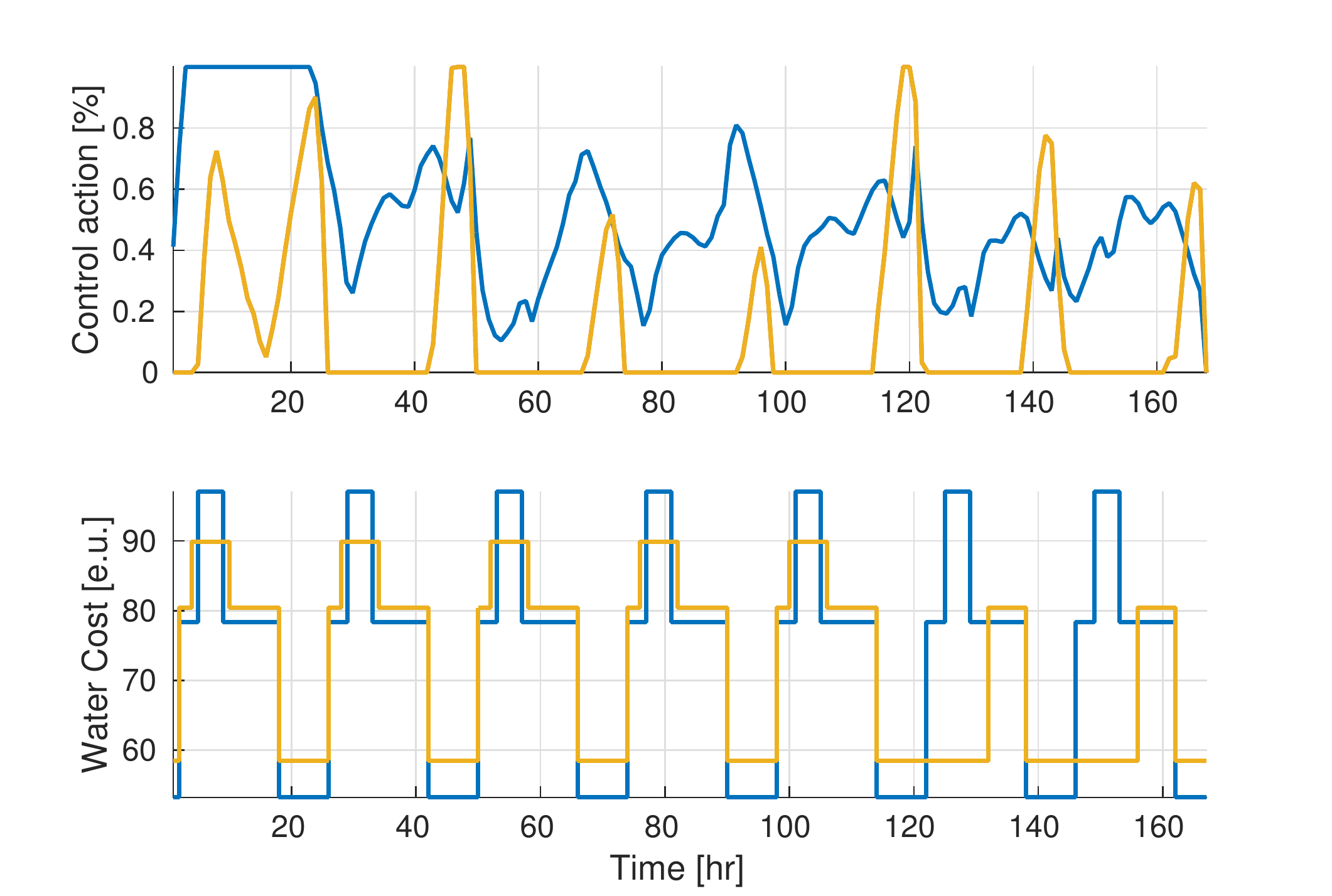}
 \caption{A pumping action using SMPC$_4$ (expressed in $\%$ of $u_{\max}$) and the 
 corresponding weighted time-varying cost $W_{\alpha}\alpha_{2,k}$ in economic units.}
 \label{fig:ControlOfController5}
\end{figure}

\subsection{Implementation details}
At every time instant $k$ we need to load onto the GPU the state measurement and 
a sequence of demand predictions (see Figure~\ref{fig:smpc_scheme}), that is ${\hat{\mathbf{d}}}_k$. 
This amounts to $8.4kB$ and is rapidly uploaded on the GPU (less than $0.034ms$).
In case we need to update the scenario-tree values, that is $\boldsymbol{\epsilon}_k$, 
and for the case of $\mathrm{SMPC}_8$ we need to upload $3.52MB$ which is done in
$3.74ms$. Therefore, the time needed to load these data on the GPU is not a limiting 
factor.

\section{Conclusions}\label{S.7}
In this paper we have presented a framework for the formulation of
a stochastic model predictive control problem for the operational 
management of drinking water networks and we have proposed a 
novel approach for the efficient numerical solution of the associated
optimization problem on a GPU.

We demonstrated the computational feasibility of the algorithm and the 
benefits for the operational management of the system in terms of 
performance (which we quantified using certain KPIs from the literature).

\appendices
\section{Elimination of input-disturbance coupling}\label{App:cost update}

In this section we discuss how the input-disturbance equality constraints
can be eliminated by a proper change of input variables and we compute the 
parameters $\beta_j^i,\hat{u}_j^i,e_j^i$ $\forall i \in \mu(j), j \in \N_{N}$ which are
then provided as input to Algorithm~\ref{algo:solve_dp}. These depend on the nominal 
demand forecasts $\hat{d}_{k+j\mid k}$ and on the time-varying
economic cost parameters $\alpha_{2,k+j}$ for $j\in\N_{[0, N-1]}$, therefore, they
need to be updated at every time instant $k$.

The affine space $\Phi_1(d)$ introduced in~\eqref{eq:Phi_1} can be written as
\begin{equation}\label{eq:new:decision}
 \Phi_1(d)=\{v \in \Re^{n_v}:u=Lv+\hat{u}(d)\},
\end{equation}
where $L \in \Re^{n_u \times n_v}$ is a \textit{full rank} matrix whose range 
spans the nullspace of $E$, i.e., for every $v \in \Re^{n_v}$, 
we have $Lv$ is in the kernel of $E$ and $\hat{u}(d)$ satisfies $E\hat{u}(d)+E_dd=0$.

Substituting $u_j^i=Lv_j^i+\hat{u}_j^i, \forall i\in\mu(j), j\in \N_{N}$ in the dynamics $\Phi_2(d)$ in~\eqref{eq:Phi_2} 
gives 
\begin{align}\label{eq:new:system:dyn} 
 \Phi_2(d)=&\{(x_{j+1}, x_j, v): x_{j+1}=Ax_j + \bar{B}v + e,\nonumber \\
 &\bar{B}=BL,e=B\hat{u}+G_d d\},
\end{align}
and we define
\begin{equation}\label{eq:appendix:dist}
 e_j^i=B\hat{u}_j^i+G_d d_j^i.
\end{equation}
Now the cost in~\eqref{eq:fz} is transformed as: 
 \begin{align}
  \sum_{j=0}^{N-1}\sum_{i=1}^{\mu(j)}  p_j^i(\ell^w(u_j^i)+\ell^{\Delta}(\Delta u_j^i))&=\nonumber\\
\sum_{j=0}^{N-1}\sum_{i=1}^{\mu(j)} p_j^i(\ell^w(v_j^i)+\ell^{\Delta}(\Delta v_j^i,\hat{u}_j^i)),\label{eq:cost:temp}
 \end{align}
 where
 \begin{subequations}
\begin{align}
\hat{R}&=W_uL,\\
  \bar{R}&=L'\hat{R},\\
  \bar{\alpha}_j&=W_{\alpha}(\alpha_1+\alpha_{2,j+k})L,\\
 \ell^w(v_j^i)&=\bar{\alpha}_j^{\prime}v_j^i,\\
 \Delta v_j^i&=v_j^i-v_{j-1}^{\parent(j,i)},\\
 \Delta \hat{u}_j^i&=\hat{u}_j^i-\hat{u}_{j-1}^{\parent(j,i)},\\
 \ell^{\Delta}(\Delta v_j^i,\Delta  \hat{u}_j^i)&=\Delta v_j^i\bar{R}\Delta v_j^i+ 2\Delta\hat{u}_j^{i \prime}\hat{R}\Delta v_j^i \label{eq:cost:delta}.
\end{align}
\end{subequations}
By substituting and expanding $\Delta v_j^i$ and $\Delta \hat{u}_j^i$ in $\ell^{\Delta}(\Delta v_j^i,\Delta \hat{u}_j^i)$ the cost 
in~\eqref{eq:cost:delta} becomes
\begin{align}\label{eq:new:stage:cost}
&\sum_{j=0}^{N-1}\sum_{i=1}^{\mu(j)}  p_j^i(\ell^w(u_j^i)+\ell^{\Delta}(\Delta u_j^i))=\nonumber\\
&\sum_{j=0}^{N-1}\sum_{i=1}^{\mu(j)} \bar{p}_j^iv_j^{i\prime}\bar{R}v_j^i-2p_j^iv_{j-1}^{\parent(j,i)\prime}\bar{R}v_j^i
         +\beta_j^{i\prime}v_j^i,
\end{align}
where 
\begin{subequations}
 \begin{align}
\bar{p}_j^i&=p_j^i+\sum_{l \in \child(j,i)}p_{j+1}^l,\\
 \beta_j^i&=p_j^i\bar{\alpha}_j+2p_j^i\hat{R} \Big(\bar{p}_j^i\hat{u}_j^i-\hat{u}_{j-1}^{\parent(j,i)}-\notag\\
  &\quad\sum_{l \in \child(j,i)}p_{j+1}^l \hat{u}_{j+1}^l\Big).\label{eq:appendix:beta}
\end{align} 
\end{subequations}
Now $\hat{u}_j^i$, $e_j^i$, $\beta_j^i $ are calculated from~\eqref{eq:new:decision},
\eqref{eq:appendix:dist} and \eqref{eq:appendix:beta} respectively. 
Using our assumption that $L$ is full-rank, we can see that $\bar{R}$ is a positive definite and symmetric matrix, 
therefore, $f$ is strongly convex.

\section{Factor step}\label{App:Factor step}
Algorithm~\ref{algo:solve_dp} solves the unconstrained minimization problem~\eqref{optim:primal:update}, that is
\begin{equation}
 z^\star=\argmin_{z}\{\langle z,H^{\prime} y\rangle +f(z)\},
\end{equation}
where $z=\{x_j^i,u_j^i\}$, $y=\{\tilde\varsigma_j^i, \tilde\zeta_j^i, \tilde\psi_j^i\}$
for $i\in \N_{[1,\mu(j)]}$ and $j\in \N_{[0,N]}$, $f(z)$ is given by~\eqref{eq:fz} and
$H$ is given by~\eqref{optim:H:operator}.
Substituting $H$ the optimization problem becomes
\begin{align}\label{optim:expanded_dual}
 z^\star&=\argmin_{z} \sum_{j=0}^{N-1}\sum_{i=1}^{\mu(j)}  
                        p_j^i(\ell^w(u_j^i)+\ell^{\Delta}(\Delta u_j^i))\nonumber\\
                      & \qquad+\tilde{\xi}_j^{i\prime}x_j^i+\tilde{\psi}_j^{i \prime}u_j^i+\delta(u_j^i|\Phi_1(d_j^i)) \notag\\
                      & \qquad+ \delta(x_{j+1}^{i}, u_j^i, x_j^{\parent(j+1,i)}|\Phi_2(d_j^i)),
\end{align}
where $\tilde{\xi}_j^i \dfn (\tilde{\varsigma}_j^i,\tilde{\zeta}_j^i)$.

The input-disturbance coupling constraints imposed by 
$\delta(u_j^i|\Phi_1(d_j^i))$ in the above problem are 
eliminated as discussed in Appendix~\ref{App:cost update}. 
This changes the input variable from $u_j^i$ to $v_j^i$ given by~\eqref{eq:new:decision} and 
the cost function as in~\eqref{eq:new:stage:cost}. 
We, therefore, replace the decision variable $z$ with $\bar{z}\dfn \{x_j^i,v_j^i\}$ and
the optimization problem~\eqref{optim:expanded_dual}
reduces to
\begin{align}\label{eq:new:fz}
 \bar{z}^\star&=\argmin_{\bar{z}} \sum_{j=0}^{N-1}\sum_{i=1}^{\mu(j)} \bar{p}_j^iv_j^{i\prime}\bar{R}v_j^i-2p_j^iv_{j-1}^{\parent(j,i)\prime}\bar{R}v_j^i\nonumber\\
                  &\qquad+\beta_j^{i\prime}v_j^i+\tilde{\xi}_{j}^{i \prime}x_{j+1}^j+\tilde{\psi}_{j}^{i \prime}Lv_j^i\nonumber\\
                  &\qquad+\delta(x_{j+1}^{i}, v_j^i, x_j^{\parent(k,i)}|\Phi_2(d_j^i)),
\end{align}
where $u_j^i=Lv_j^i+\hat{u}_j^i$. 

The above problem is an unconstrained optimization problem with quadratic stage cost which is solved using dynamic
programming~\cite{Bert00}. This method transforms the complex problem into a sequence of sub-problems solved at 
each stage. 

Using dynamic programming we find that the transformed control actions 
$v_{j}^{i\star}$ have to satisfy
\begin{align}
     v_{j}^{i\star}&=v_{j-1}^{\parent(j,i)}+
     \frac{1}{2p_{j}^i}\big (\Phi(\tilde{\xi}_j^i+q_{j+1}^i)+\Psi \tilde{\psi}_j^i\notag\\
     &+\Lambda(\beta_j^i+r_{j+1}^i)\big ),\label{eq:optimal control k}
\end{align}
where
\begin{subequations}\label{eq:optimal matrices control}
\begin{align}     
     \Lambda&= -\bar{R}^{-1},\\
     \Phi&=\Lambda\bar{B}^{\prime},\\
     \Psi&=\Lambda L.    
\end{align}
\end{subequations}
Matrix $\bar{R}$ is symmetric and positive definite, therefore, we can compute once its
Cholesky factorization so that we obviate the computation of its inverse. 

The $q_{j+1}^i$, $r_{j+1}^i$ in~\eqref{eq:optimal control k} correspond to the 
linear cost terms in the cost-to-go function at node $i$ of stage $j+1$. 
At stage $j$, these terms are updated by substituting the $v_{j}^{i\star}$ as:
\begin{subequations}
\begin{align}
   r_j^{s} &= \sum_{l \in \child(j-1,s)} \sigma_j^l  +  \bar{B}^{\prime}(\tilde{\xi}_j^l+q_{j+1}^l)+L\tilde{\psi}_j^l,\\
   q_j^{s} &= A^{\prime}\sum_{l \in \child(j-1,s)}\tilde{\xi}_j^l+q_{j+1}^l,
\end{align}\label{eq:optimal_value_function}
\end{subequations}
where $s=\parent(j,i)$.

Equations~\eqref{eq:optimal control k} and~\eqref{eq:optimal_value_function} 
form the solve step as in Algorithm~\ref{algo:solve_dp}. Matrices 
$\Lambda$, $\Phi$ and $\Psi$ are required to be computed once.

\section{Proximal operators}\label{App:Proximal}
Function $g$ in~\eqref{optim:g} is a separable sum of distance and indicator
functions and its proximal is computed according to~\eqref{eq:separable-sum-property}.
The proximal operator of the indicator of a convex closed set $C$, that is
 \[
  \chi_C(x)=\begin{cases}
        0, & \text{ if } x \in \mathcal{C}\\
        +\infty, & \text{ otherwise }
       \end{cases}
 \]
is the projection operator onto $C$, i.e., 
\begin{equation}
 \prox_{\lambda \chi_C}(v)=\proj_{\mathcal{C}}(v) = \argmin_{y\in C}\|v-y\|, 
\end{equation}

When $g$ is the distance function from a convex closed set $C$, that is
\begin{align}
 g(x)=&\mu\dist(x \mid \mathcal{C})= \inf_{y \in \mathcal{C}} \mu\| x-y \| \nonumber \\
     =& \mu \| x-\proj_{\mathcal{C}}(x)\|.\nonumber
\end{align}

Then proximal operator of $g$ given by~\cite{ComPes10}
\[
\prox_{\lambda g}{(v)}= \begin{cases}
  x+\frac{\proj_{\mathcal{C}}(x)-x}{\dist(x \mid \mathcal{C})}, & \text{ if } \dist(x\mid \mathcal{C})>\lambda\mu\\
  \proj_{\mathcal{C}}(x), & \text{ otherwise }
 \end{cases}
 \]

\section*{Acknowledgment}
This work was financially supported by the EU FP7
research project EFFINET ``Efficient Integrated Real-time 
monitoring and Control of Drinking Water Networks,''
grant agreement no. 318556.

\ifCLASSOPTIONcaptionsoff
  \newpage
\fi

\bibliographystyle{ieeetr}
\bibliography{waterlit}

\end{document}